\documentclass[11pt,letterpaper,reqno]{amsart}
\usepackage{mathrsfs}
\usepackage{cite}
\usepackage{graphicx}
\usepackage{enumitem}
\usepackage{microtype}
\usepackage[colorlinks, linkcolor=blue,citecolor=blue]{hyperref}

\numberwithin{equation}{section}
\newcommand{\R}{\mathbb R}

\newcommand{\Q}{\mathbb Q}
\newcommand{\N}{\mathbb N}

\newcommand{\norm}[1]{\left\lVert #1\right\rVert}
\newcommand{\abs}[1]{\left\lvert #1\right\rvert}

\newtheorem{theo}{Theorem}[section]
\newtheorem{lem}[theo]{Lemma}
\newtheorem{prop}[theo]{Proposition}

\theoremstyle{definition}

\newtheorem{rem}[theo]{Remark}

\numberwithin{equation}{section}

\begin{document}

\title[Inhomogeneous Fiber Refinement of Gallagher's Theorem]{Inhomogeneous Fiber Refinement of Gallagher's Theorem}

\author[Q. Li]{Qian Li}
\address{School of Mathematics and Statistics, Wuhan University of Technology}
\email{l.qian@whut.edu.cn}

\author[B. Tan]{Bo Tan}
\address{School of Mathematics and Statistics, Huazhong University of Science and Technology}
\email{tanbo@hust.edu.cn}

\author[Q.-L. Zhou$^{\ast}$]{Qing-Long Zhou$^{\ast}$}
\address{School of Mathematics and Statistics, Wuhan University of Technology}
\email{zhouql@whut.edu.cn}

\thanks{*Corresponding author.}
\thanks{This work is supported by National Key R$\&$D Program of China (No. 2024YFA1013700) and NSFC (Nos. 12671113, 12331005 and 12601156), and the Fundamental Research Funds for the Central Universities.}

\begin{abstract}
Let \(s>0\), \(\alpha,\gamma\in\mathbb{R}\), and let \(\psi\colon \mathbb{N}\to\mathbb{R}_{>0}\) be a non-increasing function. Consider the limsup set
\[
    W_s(\psi,\alpha,\gamma) = \left\{ 
        \beta\in[0,1] \colon 
        \norm{n\beta-\gamma} < \frac{\psi(n)}{\norm{n\alpha}^{\,s}} 
        \ \text{ for infinitely many } n\in\mathbb{N}
    \right\}.
\]
Chow [Duke J. Math. 2018] settled the measure-theoretic result for \(s=1\). In this paper, we investigate the Lebesgue measure of \(W_s(\psi,\alpha,\gamma)\) for general \(s>0\).

\end{abstract}

\maketitle

\section{Introduction}

\subsection{Metrical Diophantine approximation}

A central theme in Diophantine approximation is to understand how well real numbers can be approximated by rationals. At the heart of metric number theory lies the study of the size---quantified via Lebesgue measure, Hausdorff dimension, Fourier dimension, etc.---of Diophantine sets defined by prescribed approximation conditions.

Given a function \(\psi\colon \mathbb{N} \to [0,\infty)\) and a number \(\alpha \in [0,1]\) (typically irrational), Diophantine approximation studies the existence of infinitely many pairs \((m,n) \in \mathbb{Z} \times \mathbb{N}\) satisfying
\begin{equation}\label{eq:diophantine}
    \left| \alpha - \frac{m}{n} \right| < \frac{\psi(n)}{n}.
\end{equation}
For the critical case \(\psi(n) = 1/n\), Dirichlet's Theorem guarantees that \eqref{eq:diophantine} holds for all \(\alpha \in [0,1]\). While this bound is asymptotically sharp for badly approximable numbers, Khintchine's Theorem reveals that we can improve upon it for Lebesgue almost every \(\alpha \in [0,1]\).

Define the set of well-approximable numbers as
\[
    W(\psi,\delta) := \left\{ \alpha \in [0,1]\colon \norm{n\alpha-\delta} < \psi(n) \text{ for i.m. } n \in \mathbb{N} \right\}.
\]
Here and in what follows ``i.m.'' stands for ``infinitely many'', and $\|\alpha\|$ denotes the distance from $\alpha$ to the nearest integer.

\begin{theo}[Khintchine, \cite{Khintchine1924}]
Let \(\psi\) be a monotonically decreasing function. Then the Lebesgue measure \(\mathcal{L}(W(\psi,0))\) is given by
\[
    \mathcal{L}\bigl(W(\psi,0)\bigr) = 
    \begin{cases} 
        0, & \text{if } \displaystyle\sum_{n \in \mathbb{N}} \psi(n) < \infty, \\[1.4em]
        1, & \text{if } \displaystyle\sum_{n \in \mathbb{N}} \psi(n) = \infty.
    \end{cases}
\]
\end{theo}

Duffin and Schaeffer \cite{Duffin-Schaeffer1941} later demonstrated that the monotonicity condition is essential: they constructed a \(\psi\) supported on highly composite integers such that
\[
    \sum_{n=1}^{\infty} \psi(n) = \infty, \qquad \mathcal{L}\bigl(W(\psi,0)\bigr) = 0.
\]
They conjectured that for almost all \(\alpha\), the inequality \(|n\alpha - m| < \psi(n)\) has infinitely many coprime solutions \((n, m)\) if and only if
\begin{equation}\label{eq:phi-psi-divergence}
    \sum_{n=1}^{\infty} \frac{\phi(n)}{n} \psi(n) = \infty,
\end{equation}
where \(\phi\) is Euler's totient function. Following contributions by Gallagher \cite{G61}, Erd\H{o}s\cite{E70}, Vaaler\cite{V78}, and Pollington--Vaughan\cite{PV90}, the conjecture was resolved affirmatively by Koukoulopoulos and Maynard \cite{Koukoulopoulos-Maynard2020}.

The inhomogeneous analogue of Khintchine's theorem was established by Sz\"usz \cite{Szusz1958}. 
\begin{theo}[Sz\"usz, \cite{Szusz1958}]\label{thm:inhomogeneous-Khintchine}
Let \(\psi\colon \mathbb{N}\to\mathbb{R}_{+}\) be monotonically decreasing, and let \(\delta\in\mathbb{R}\) be fixed. Then
\[
    \mathcal{L}\bigl(W(\psi,\delta)\bigr)
    =
    \begin{cases}
        0, & \text{if } \displaystyle\sum_{n\in\mathbb{N}} \psi(n) < \infty, \\[1.4em]
        1, & \text{if } \displaystyle\sum_{n\in\mathbb{N}} \psi(n) = \infty.
    \end{cases}
\]
\end{theo}
The non-monotonic case is more intricate: counterexamples to relaxing monotonicity were given in \cite{CHPR25,R17}. Recently, 
The coprime-numerator assertion fails in general. Hauke-Treuer, Maynard and Pollington \cite{HMP26} constructed counterexamples for all nonzero rational shifts and certain Liouville shifts; independently, He and Liao \cite{HL26} obtained counterexamples for all nonzero rational shifts and a residual set of real shifts.

An open problem in this direction is the so-called \emph{weak inhomogeneous Duffin--Schaeffer conjecture} (Conjecture 1.22, \cite{Chow-Technau2024}). It asserts that if the divergence condition \eqref{eq:phi-psi-divergence}
holds, then for almost all \(\alpha \in [0,1]\), the inequality \(\norm{n\alpha - \delta} < \psi(n)\) should have infinitely many solutions \(n \in \mathbb{N}\). 
Partial progress has been made in various aspects. For instance, it was solved for a wide class of non-monotonic functions in \cite{Chow2018,CT19,Chow-Technau2024}. Results with extra divergence were established by Yu \cite{Yu2019,Yu2021}. 
More recently, OpenAI has announced a proof of the weak inhomogeneous Duffin–Schaeffer conjecture \cite{openai2026}.

\begin{theo}[OpenAI, \cite{openai2026}]\label{thm:OpenAI}
Let \(\delta \in \mathbb{R}\), and let \(\psi\colon \mathbb{N} \to [0,\infty)\). Then 
\[
    \mathcal{L}\bigl(W(\psi,\delta)\bigr)=1 \quad \text{if} \quad \sum_{n=1}^{\infty} \frac{\phi(n)}{n}\psi(n) = \infty.
\]
\end{theo}

\begin{rem}
While the divergence case is now fully resolved, a complete zero-one law for \(W(\psi, \delta)\) remains unknown. 
\end{rem}

\subsection{Multiplicative Diophantine approximation}
Multiplicative Diophantine approximation is closely tied to Littlewood's conjecture, which asserts that
\begin{equation*}
\liminf_{n\to\infty} n\|n\alpha\|\|n\beta\|=0
\end{equation*}
for all \(\alpha,\beta\in[0,1]\). 
Einsiedler, Katok, and Lindenstrauss \cite{EKL06} proved that the exceptional set has Hausdorff dimension zero. However, the conjecture remains open, even for natural cases like \((\alpha, \beta) = (\sqrt{2}, \sqrt{3})\). To date, the only non-trivial resolved instances are those where \(\alpha\) and \(\beta\) lie in the same cubic field \cite{CS55,P61}.

Turning to the metric theory, we define the set of multiplicative well-approximable numbers as
\[
    W^{\times}(\psi) := \left\{ (\alpha,\beta) \in [0,1]\times[0,1]\colon \|n\alpha\| \|n\beta\|< \psi(n) \text{ for i.m. } n \in \mathbb{N} \right\}.
\]
\begin{theo}[Gallagher, \cite{Gallagher1962}]
Let \(\psi\colon \mathbb{N}\to\mathbb{R}_{+}\) be monotonically decreasing. Then
\[
    \mathcal{L}\bigl(W^{\times}(\psi)\bigr)
    =
    \begin{cases}
        0, & \text{if } \displaystyle\sum_{n\in\mathbb{N}} \psi(n)\log n < \infty, \\[1.4em]
        1, & \text{if } \displaystyle\sum_{n\in\mathbb{N}} \psi(n)\log n = \infty.
    \end{cases}
\]
\end{theo}
Without assuming monotonicity, Beresnevich, Haynes and Velani \cite{BHV13} proved a dichotomy for the Lebesgue measure of \(W^{\times}(\psi)\) under additional assumptions. By removing these conditions, Frühwirth and Hauke \cite{FH24} established the following: for almost all \((\alpha, \beta) \in [0, 1]^2\), there exist infinitely many \(n\) such that 
\(|n\alpha-m_1| |n\beta-m_2|< \psi(n)\) with \(m_1, m_2\) coprime to \(n\), if and only if the series \(\sum_n \frac{\phi(n)\psi(n)}{n} \log \left( \frac{n}{\phi(n)\psi(n)} \right)\) diverges.

Inhomogeneous multiplicative Diophantine approximation yields a natural yet nontrivial generalization of this theory. Fix \(\delta, \gamma \in \mathbb{R}\),
Shapira \cite{Shapira2011} proved that for almost every pair $(\alpha,\beta)$,
\[
\liminf_{n\to\infty} n\|n\alpha-\delta\|\|n\beta-\gamma\| = 0.
\]
 This was subsequently improved by Gorodnik and Vishe \cite{Gorodnik-Vishe2018}, who showed that there exists \(\varepsilon > 0\) such that almost all \((\alpha, \beta) \in [0,1]^2\) satisfy
\[
    \liminf_{n\to\infty} (\log_5 n)^\varepsilon n\|n\alpha - \gamma\| \cdot \|n\beta - \delta\| = 0,
\]
where \(\log_5\) denotes the fifth iterate of the function \(x \mapsto \max\{1, \log x\}\). 
While the above results support Littlewood's conjecture, Adiceam, Nesharim, and Lunnon \cite{ANL21} provided a counter-indication by proving that its function field analogue is false.

Fibre results are especially relevant to Littlewood-type problems
because almost-everywhere results in the whole plane can avoid the
most arithmetic pairs. Recall that the set of badly approximable numbers is defined as
\[
    \mathbf{Bad} := \left\{ \alpha \in [0,1] \colon \liminf_{n \to \infty} n\|n\alpha\| > 0 \right\}.
\]
Borel--Cantelli lemma yields that $\mathbf{Bad}$ has Lebesgue measure zero, while by the Jarník--Besicovitch theorem, it has full Hausdorff dimension. 
Pollington and Velani \cite{Pollington-Velani2000} proved that for every \(\alpha\in\mathbf{Bad}\), there exists a set \(G\subseteq\mathbf{Bad}\) with \(\dim_{\mathrm H}G=1\) such that
\[
    \liminf_{n\to\infty} n\log n\,\|n\alpha\|\,\|n\beta\|\leq1
    \qquad\text{for every }\beta\in G.
\]
Chow and Zafeiropoulos \cite{Chow-Zafeiropoulos2021} extended this result to the fully inhomogeneous setting and relaxed the assumption on \(\alpha\). Related developments include the
Erd\H{o}s--Vaaler-type fibre counting results of Yu \cite{Yu2023}
and the almost-everywhere multiplicative results on planar lines
obtained by Chow and Yang \cite{Chow-Yang2024}.

Here we recall that an irrational real number \(\alpha\) is \textit{Liouville} if
\[
    \liminf_{n\to\infty} n^w \|n\alpha\| = 0
\]
for all \(w > 0\). 
The  exponent of approximation of \(\alpha \in [0,1]\) is defined as
\begin{equation*}
    w(\alpha) := \sup\{w > 0\colon  \|n\alpha\| < n^{-w} \text{ for i.m. } n \in \mathbb{N}\}.
\end{equation*}
Note that, for an irrational \(\alpha\), by definition, \(w(\alpha) < \infty\) if and only if \(\alpha\) is not Liouville.  Beresnevich, Haynes and Velani \cite{BHV13} proved that if \(\alpha\) is irrational and non-Liouville,  and assuming the Duffin--Schaeffer conjecture, then
\[
    \liminf_{n\to\infty} n(\log n)^2 \|n\alpha - \delta\| \cdot \|n\beta\| = 0
\]
for almost all \(\beta \in [0,1]\). 
Later, Chow \cite{Chow2018} proved this result unconditionally (without assuming the Duffin--Schaeffer conjecture) and in a stronger, function-based form. By finding large generalised arithmetic progressions within inhomogeneous Bohr sets, he showed that if \(\alpha\) is irrational and non-Liouville, and \(\psi\) is nonincreasing, then
\[
    \sum_{n=1}^{\infty}\psi(n)\log n=\infty
    \quad\Longrightarrow\quad
    \|n\alpha-\delta\|\,\|n\beta\|<\psi(n)
    \ \text{for i.m. }n,
\]
for almost every \(\beta\). Furthermore, by combining this with Corollary 2.1 in \cite{BHV13}, Chow established a complete zero-one law$\colon$ if \(\alpha\in \mathbf{Bad}\) and \(\psi\colon \mathbb{N} \to \mathbb{R}_{\ge 0}\) is a nonincreasing function, then the Lebesgue measure of the set
\[
    \left\{\beta \in [0,1] :   \|n\beta\| < \frac{\psi(n)}{\|n\alpha\|} \text{ for i.m. } n \in \mathbb{N}\right\}
\]
is 0 or 1  according to the convergence or divergence of the series \(\sum_{n}\psi(n) \log n\). The fully inhomogeneous extension was established in the recent work of Chow and Technau \cite{Chow-Technau2024}. In particular, building upon Corollary 2.1 in \cite{BHV13}, they  proved that if \(\alpha\in \mathbf{Bad}\),  and \(\psi\) is nonincreasing, then the Lebesgue measure of the set
\[
    \left\{\beta \in [0,1] : \|n\beta-\gamma\| < \frac{\psi(n)}{\|n\alpha\|} \text{ for i.m. } n \in \mathbb{N}\right\}
\]
is \(0\) or \(1\) according to the convergence or divergence of the series \(\sum_{n}\psi(n) \log n\).

In this paper, we further extend the inhomogeneous approximation problem by incorporating a weight into the approximation condition. More precisely, for an irrational number \(\alpha\), a shift \(\gamma \in \mathbb{R}\), a weight parameter \(s>0\), and a nonincreasing function \(\psi\colon\mathbb{N}\to\mathbb{R}_{>0}\), we study the Lebesgue measure of the set 
\begin{equation*}
 W_s(\psi,\alpha,\gamma)
 :=\left\{\beta\in[0,1]:
 \norm{n\beta-\gamma}<\frac{\psi(n)}{\norm{n\alpha}^{s}}
 \text{ for i.m. } n \in \mathbb{N}\right\}.
\end{equation*}
Our main results read:
\begin{theo}\label{thm:metric-subcritical}
Let $0<s<1$, $\alpha\in[0,1]\setminus\Q$ and $\gamma\in\R$. 
Let \(\psi\colon \mathbb{N}\to\mathbb{R}_{+}\) be a nonincreasing function. Then
\[
    \mathcal{L}\bigl(W_s(\psi,\alpha,\gamma)\bigr)
    =
    \begin{cases}
        0, & \text{if } \displaystyle\sum_{n\in\mathbb{N}} \psi(n)\log n < \infty\text{ and }w(\alpha)<1/s, \\[1.4em]
        1, & \text{if } \displaystyle\sum_{n\in\mathbb{N}} \psi(n)\log n = \infty .
    \end{cases}
\]
\end{theo}


\begin{theo}\label{thm:metric-supercritical-initial}
Let \(s>1\), \(\alpha\in[0,1]\setminus\mathbb{Q}\) and \(\gamma\in\mathbb{R}\). 
Let \(\psi\colon \mathbb{N}\to\mathbb{R}_{+}\) be a nonincreasing function. 
\begin{enumerate}\setlength{\itemsep}{0.6em} 
\item If for some \(\varepsilon_0>0\),
\[
    \sum_{n=1}^{\infty}\psi(n)n^{w(\alpha)s+\varepsilon_0-1}<\infty,
\]
then \(\mathcal{L}\bigl(W_s(\psi,\alpha,\gamma)\bigr)=0\).

\item  If
\[
    \sum_{n=1}^{\infty}\psi(n)\frac{n^{s-1}}{\log\log n}=\infty,
\]
then \(\mathcal{L}\bigl(W_s(\psi,\alpha,\gamma)\bigr)=1\).
\end{enumerate}
\end{theo}

\medskip

Let us make the following remarks regarding Theorems \ref{thm:metric-subcritical} and \ref{thm:metric-supercritical-initial}:
\begin{itemize}
    \item Khintchine's Theorem and Dirichlet's Theorem yields that \(\mathcal{L}\bigl(\{\alpha \in [0,1]\colon w(\alpha) = 1\}\bigr) = 1\). 
If \(s < 1\), the condition \(w(\alpha) = 1 < 1/s\) is automatically satisfied for almost every \(\alpha\). 

\item Theorem \ref{thm:metric-supercritical-initial} does not yield a complete zero-one law. The parameter \(\varepsilon_0\)  originates from the limit superior in the definition of \(w(\alpha)\). If   \(q_{k+1} \ll q_{k}^{w(\alpha)}\) holds eventually, this parameter can be eliminated. Namely, if \(\sum \psi(n) n^{w(\alpha)s-1} < \infty\), then \(\mathcal{L}(W(\psi, \alpha, \gamma)) = 0\).
The factor \(1/\log\log n\) in the divergence condition arises naturally from Landau's classical estimate for Euler's totient function\cite{Landau1909}: \(\phi(n)/n \gg 1/\log\log n\). This bound is optimal, as the limit inferior of \(\phi(n)/(n\log\log n)\) is the positive constant \(e^{-c}\) (where \(c \approx 0.5772\) is the Euler--Mascheroni constant).  

\end{itemize}

\medskip

\noindent\textbf{Application.} For \(a, b > 0\) and a fixed irrational \(\alpha\), we define the weighted multiplicatively well-approximable fiber as
\[
    W^{\times}_{a,b}(\psi, \alpha, \gamma) := \left\{ \beta \in [0,1] \colon \norm{n\alpha}^a \norm{n\beta - \gamma}^b < \psi(n) \text{ for i.m. } n \in \mathbb{N} \right\}.
\]
By a direct substitution, this set is equivalent to the weighted inhomogeneous set:
\[
    W^{\times}_{a,b}(\psi, \alpha, \gamma)= W_{a/b}(\psi^{1/b}, \alpha, \gamma).
\]

\medskip

\noindent\textbf{Organization.} The paper is organized as follows. In Section 2, we collect some known auxiliary results that will be used in our proofs. In Section 3, we estimate the reciprocal sums \(\sum_{1\le n\leq N}\|n\alpha\|^{-s}\) for \(0<s<1\) and \(s>1\), which serve as a key ingredient in our measure estimates. Sections 4 and 5 are devoted to the proofs of Theorems \ref{thm:metric-subcritical} and \ref{thm:metric-supercritical-initial}, respectively.

\medskip

\noindent\textbf{Notation.} We use Vinogradov's notation: \(A\ll B\) means that \(A\le CB\) for some constant \(C>0\). We write \(A\gg B\) if \(B\ll A\), and \(A\asymp B\) if \(A\ll B\ll A\).

\section{Preliminaries}

In this section, we present some auxiliary results required for the proofs of our main theorems. 


Let $\alpha\in[0,1]\setminus\mathbb{Q}$, and write
\[
\alpha=[a_1,a_2,\ldots]
\]
for its continued-fraction expansion. Further, we write the finite truncations
\[
    [a_1, a_2, \dots, a_j] := \frac{p_j}{q_j}
\]
to be the \(n\)-th convergent to \(\alpha\). The following relations can be found in \cite{RS1992}:
\begin{equation*}
q_{j+1}=a_{j+1}q_j+q_{j-1},
\qquad
p_{j+1}=a_{j+1}p_j+p_{j-1},
\end{equation*}
and
\begin{equation*}
\frac{1}{2q_{j+1}}
<
\abs{q_j\alpha-p_j}
<
\frac{1}{q_{j+1}}.
\end{equation*}
We will often make use of the quantities
\begin{equation*}
D_k = q_k \alpha - p_k \qquad (k \ge 0).
\end{equation*}
Then
\begin{equation*}
D_0 = \{\alpha\},\qquad D_k = (-1)^k \|q_k \alpha\| \qquad (k\ge 1),
\end{equation*}
and
\begin{equation*}
\frac12 \le q_{k+1}|D_k| \le 1 \qquad (k\ge 0).
\end{equation*}

The following lemma provides a useful criterion for computing the Diophantine exponent \(\omega(\alpha)\) via the partial quotients of its continued fraction expansion.

\begin{lem}[\hspace{-0.03em}\cite{Harman1998}]\label{lem:omega-def}
Let \(\alpha \in [0,1]\setminus\mathbb{Q}\), and let \(\{q_k\}_k\) be the sequence of its continued fraction denominators. Then
\[
    \omega(\alpha) = \limsup_{k\to\infty} \frac{\log q_{k+1}}{\log q_k}.
\]
\end{lem}

The following lemma presents the \textit{Ostrowski expansion} (or \textit{Ostrowski numeration}) of the integers.

\begin{lem}[\hspace{-0.03em}\cite{Ostrowski1922}]\label{lem:Ostrowski}
Let \(\alpha=[a_1,a_2,\ldots] \in [0,1]\setminus\mathbb{Q}\). Then, for every \(n\in\mathbb{N}\), there is a unique integer \(K\ge 0\) such that
\[
    q_K \le n < q_{K+1},
\]
and a unique sequence \(\{c_{k+1}\}_{k=0}^{\infty}\) of integers such that
\[
    n = \sum_{k=0}^{\infty} c_{k+1} q_k,
\]
where
\begin{align*}
    0 \le c_1 < a_1 \quad &\text{and} \quad 0 \le c_{k+1} \le a_{k+1} \quad (\forall\,k\ge 1), \\
    &c_k = 0 \quad \text{whenever} \quad c_{k+1}=a_{k+1} \text{ with } k\ge 1, \\
    &c_{k+1}=0 \quad (\forall\,k>K).
\end{align*}
\end{lem}

Using the unique Ostrowski expansion \(n = \sum_{k=0}^{\infty} c_{k+1} q_k\) provided by Lemma~\ref{lem:Ostrowski}, let \(m \ge 0\) be the minimal integer such that \(c_{m+1} \neq 0\). The following lemma expresses \(\|n\alpha\|\) in terms of the corresponding continued-fraction errors \(D_k\).

\begin{lem}[\hspace{-0.03em}\cite{Beresnevich-Haynes-Velani2020}, Lemma 4.1]\label{lem:homogeneous-norm-estimate}
Let $\alpha\in[0,1]\setminus\mathbb{Q}$, $n\in\mathbb{N}$, and, with reference to Lemma \ref{lem:Ostrowski}, let $m$ be the smallest integer such that $c_{m+1}\neq 0$.
If $m\geq 2$ then
\begin{equation}\label{eq:homogeneous-norm}
\| n \alpha\| =\left|\sum_{k=m}^{\infty} c_{k+1} D_{k}\right|=\operatorname{sgn}\left(D_{m}\right)\cdot \sum_{k=m}^{\infty} c_{k+1} D_{k}.
\end{equation}
Also if $m=1$ and $\{\alpha\}<1/2$, then we have \eqref{eq:homogeneous-norm}.
In all other cases, $\|n\alpha\| \geq |D_{2}|$, which is a positive constant depending on $\alpha$ only.
\end{lem}

Lemma \ref{lem:homogeneous-norm-estimate} underpins the following two-sided estimate for $\|n\alpha\|$.

\begin{lem}[\hspace{-0.03em}\cite{Beresnevich-Haynes-Velani2020}, Lemma 4.2]\label{lem:two-sided-norm-bound}
Let $\alpha\in[0,1)\setminus\mathbb{Q}$ and $n\in\mathbb{N}$, and let $m$ be the smallest integer such that $c_{m+1}\neq 0$. If $\|n\alpha\| < |D_{2}|$ then $m\geq 2$, and if $m\geq 2$ then
\begin{equation*}
(c_{m+1}-1)|D_{m}|+(a_{m+2}-c_{m+2})|D_{m+1}|
\leq \|n\alpha\|
\leq (c_{m+1}+1)|D_{m}|.
\end{equation*}
\end{lem}

Fix \(\ell \ge 0\). Let \(d_1, \dots, d_{\ell+1}\) be non-negative integers satisfying
\[
    d_1 < a_1, \qquad d_{k+1} \le a_{k+1} \quad (1 \le k \le \ell),
\]
and
\[
    d_k = 0 \quad \text{whenever } d_{k+1} = a_{k+1}, \quad 1 \le k \le \ell.
\]
We define \(A(d_1, \dots, d_{\ell+1})\) to be the set of all integers \(n \in \mathbb{N}\) whose Ostrowski expansion begins with the fixed coefficients \(d_1, \dots, d_{\ell+1}\). 
The corresponding truncated set is defined as
\begin{equation}\label{cylinder}
     A_N(d_1,\dots,d_{\ell+1})
    := A(d_1,\dots,d_{\ell+1}) \cap \{1, \dots, N\}.
\end{equation}

The following lemma provides a precise bound on the cardinality of these truncated sets.

\begin{lem}[\hspace{-0.03em}\cite{Beresnevich-Haynes-Velani2020}, Lemma 5.3]\label{lem:card-bound-Ostrowski-A_N}
Let \(N \ge 1\) be an integer. Provided that \(A_N(d_1,\dots,d_{\ell+1})\) is non-empty, i.e., \(\sharp A_N(d_1,\dots,d_{\ell+1}) \ge 1\), we have
\begin{equation*}
    \frac{N}{3q_{\ell+1}} \le \sharp A_N(d_1,\dots,d_{\ell+1}) \le \frac{3N}{q_{\ell+1}} + 1.
\end{equation*}
\end{lem}



\section{Estimates for reciprocal sums involving $\|n\alpha\|^{-s}$}

\medskip

This section is devoted to estimating the reciprocal sum
\[
    \sum_{1\le n\leq N}\frac{1}{\|n\alpha\|^s}.
\]
The core idea originates from Beresnevich, Haynes and Velani \cite{Beresnevich-Haynes-Velani2020}. For a given \(N\), let \(K\) be the integer uniquely determined by
\[
    q_K \le N < q_{K+1}.
\]
Throughout this section, we assume \(N \ge q_3\). This choice of \(K\) allows us to partition the integers \(1\leq n\leq N\) into two disjoint sets according to their residues modulo \(q_K\): the \textit{regular} class \(\Omega_{N,1}\) and the \textit{exceptional} class \(\Omega_{N,2}\), defined as
\begin{align*}
    \Omega_{N,1} &:= \left\{ n \in \mathbb{N} : 1 \le n \le N, \ n \not\equiv 0 \pmod{q_K} \text{ and } n \not\equiv q_{K-1} \pmod{q_K} \right\}, \\
    \Omega_{N,2} &:= \left\{ n \in \mathbb{N} : 1 \le n \le N, \ n \equiv 0 \pmod{q_K} \text{ or } n \equiv q_{K-1} \pmod{q_K} \right\}.
\end{align*}
The exceptional class \(\Omega_{N,2}\) requires special attention, since it corresponds to particularly small values of \(\|n\alpha\|\).

\begin{prop}\label{pro:s<1}
Let \(0<s<1\), \(\alpha\in\mathbb{R}\setminus\mathbb{Q}\) and \(N\ge q_3\). Let \(K\) be the largest integer satisfying \(q_K\le N\). Then
\begin{equation*}
    \sum_{n \in \Omega_{N,1}} \frac{1}{\norm{n\alpha}^{s}} \asymp N 
    \qquad \text{and} \qquad 
    \sum_{n \in \Omega_{N,2}} \frac{1}{\norm{n\alpha}^{s}} \asymp \frac{N^{1-s}q_{K+1}^{s}}{q_K^{1-s}}.
\end{equation*}
\end{prop}

\begin{prop}\label{pro:s>1}
Let \(s>1\), \(\alpha\in\mathbb{R}\setminus\mathbb{Q}\) and \(N\ge q_3\). Let \(K\) be the largest integer satisfying \(q_K\le N\). Then
\begin{equation*}
    \sum_{n \in \Omega_{N,1}} \frac{1}{\norm{n\alpha}^{s}} \asymp Nq_K^{s-1}
    \qquad \text{and} \qquad
    \sum_{n \in \Omega_{N,2}} \frac{1}{\norm{n\alpha}^{s}} \asymp q_{K+1}^{s}.
\end{equation*}
\end{prop}

The proofs of the two propositions are structured according to the partition of the integers into regular and exceptional classes. Specifically, Subsection \ref{subsec:regular} is devoted to the estimates for the regular class \(\Omega_{N,1}\), while the contribution from the exceptional class \(\Omega_{N,2}\) will be treated in Subsection \ref{subsec:exceptional}.

\subsection{Estimates for the regular  class}\label{subsec:regular}

We begin by estimating the sum over the regular class \(\Omega_{N,1}\).  We split the sum into two parts: \(\norm{n\alpha} < |D_2|\) and \(\norm{n\alpha} \ge |D_2|\):
\begin{align*}
    \sum_{n \in \Omega_{N,1}} \frac{1}{\norm{n\alpha}^{s}}
    =
    \sum_{\substack{n \in \Omega_{N,1} \\ \norm{n\alpha} < |D_2|}} \frac{1}{\norm{n\alpha}^{s}}
    +
    \sum_{\substack{n \in \Omega_{N,1} \\ \norm{n\alpha} \ge |D_2|}} \frac{1}{\norm{n\alpha}^{s}}.
\end{align*}

\noindent\underline{\textbf{Upper bounds for \(\sum_{n \in \Omega_{N,1}}\norm{n\alpha}^{-s}\).}}
For integers \(u\) and \(a\) satisfying \(1\leq u\leq K+1\) and \(1\leq a\leq a_{u+1}\), define
\begin{equation*}
    C_N(u,a) := \left\{ 1\le n\le N \colon 
    \begin{aligned}
        &\text{the Ostrowski coefficients } \{c_{k+1}\}_{k=0}^{\infty} \text{ of } n \text{ satisfy} \\
        &c_1=\cdots=c_{u-1}=0, \text{ and either} \\
        &\text{(i) } c_u=0 \text{ and } c_{u+1}=a+1, \text{ or} \\
        &\text{(ii) } c_u=1 \text{ and } c_{u+1}=a_{u+1}-a
    \end{aligned}
    \right\}.
\end{equation*}

Clearly, each \(C_N(u,a)\) is a subset of \(\{1,\dots,N\}\) by definition. Furthermore, the uniqueness of the Ostrowski expansion immediately implies that these sets are pairwise disjoint; that is, for any two distinct pairs \((u_1,a) \neq (u_2,a')\), we have \(C_N(u_1,a) \cap C_N(u_2,a') = \emptyset\). We now focus on the integers \(n\) for which \(\norm{n\alpha}<|D_2|\). As the following lemma shows, these integers are fully covered by this disjoint union.

\begin{lem}\label{lem:cover}
\begin{equation*}
 \left\{1\le n\le N:\norm{n\alpha}<|D_2|\right\}
 \subseteq\bigsqcup_{u=2}^{K+1}\bigsqcup_{a=1}^{a_{u+1}} C_N(u,a).
\end{equation*}
\end{lem}
\begin{proof}
Take \(n\) in the set on the left, and let \(m\) be the smallest index such that \(c_{m+1}\neq 0\) in its Ostrowski expansion. By Lemma \ref{lem:two-sided-norm-bound}, we have \(m\ge 2\). Since \(n\le N<q_{K+1}\), we also have \(m\le K\). Thus, \(2\le m\le K\).

We now distinguish two cases based on the value of \(c_{m+1}\):

\begin{itemize}
    \item \textbf{Case 1: \(c_{m+1} > 1\).} 
    We set \(u=m\) and \(a=c_{m+1}-1\). Then \(1 \le a \le a_{m+1}\), and it follows directly from the definition that \(n\in C_N(u,a)\).
    
    \item \textbf{Case 2: \(c_{m+1} = 1\).} 
    We set \(u=m+1\) and \(a=a_{m+2}-c_{m+2}\). Since \(0 \le c_{m+2} < a_{m+2}\), we have \(1 \le a \le a_{m+2}\). Hence, \(n\in C_N(m+1,a)\).
\end{itemize}

In either case, we obtain \(2\le u\le K+1\) and \(1\le a\le a_{u+1}\). This proves that \(n\) belongs to the disjoint union on the right.
\end{proof}

By Lemma \ref{lem:cover}, any integer \(n\) satisfying \(\norm{n\alpha}<|D_2|\) must belong to some \(C_N(u,a)\) with \(2\le u\le K+1\) and \(1\le a\le a_{u+1}\).  In particular, if \(n \in \Omega_{N,1}\), the range of \(u\) can be further restricted. The following lemma shows that the boundary cases \(u \in \{K, K+1\}\) always fall into the exceptional class \(\Omega_{N,2}\). Consequently, for \(n \in \Omega_{N,1}\) with \(\norm{n\alpha}<|D_2|\), the associated parameters must satisfy \(2 \le u \le K-1\) and \(1 \le a \le a_{u+1}\).

\begin{lem}\label{lem:terminal}
If \(n \in \bigsqcup_{u=K}^{K+1} \bigsqcup_{a=1}^{a_{u+1}} C_N(u,a)\), then
\begin{equation*}
    n \equiv 0 \pmod{q_K} \quad \text{or} \quad n \equiv q_{K-1} \pmod{q_K}.
\end{equation*}
\end{lem}

\begin{proof}
If \(n \in C_N(K,a)\), then \(n \le N < q_{K+1}\), and its Ostrowski expansion can only be \((a+1)q_K\) or \(q_{K-1} + (a_{K+1}-a)q_K\). Thus, \(n \equiv 0 \pmod{q_K}\) or \(n \equiv q_{K-1} \pmod{q_K}\), respectively.

If \(n \in C_N(K+1,a)\), its Ostrowski expansion can only be \(n = q_K + (a_{K+2}-a)q_{K+1}\). Since \(n \le N < q_{K+1}\), this forces \(n=q_K\), giving \(n \equiv 0 \pmod{q_K}\). 

\end{proof}

\begin{lem}\label{R_1u}
  We have
  \[
     \sum_{\substack{n \in \Omega_{N,1} \\ \norm{n\alpha} < |D_2|}} \frac{1}{\norm{n\alpha}^{s}}\ll  
    \begin{cases}
        N, & \text{if } 0<s<1, \\[1.4em]
        Nq_{K}^{s-1}, & \text{if } s>1.
    \end{cases}
\]
\end{lem}

\begin{proof}

First, by Lemmas \ref{lem:cover} and \ref{lem:terminal}, we obtain
\[ \sum_{\substack{n \in \Omega_{N,1} \\ \norm{n\alpha} < |D_2|}} \frac{1}{\norm{n\alpha}^{s}}\le \sum_{u=2}^{K-1}\sum_{a=1}^{a_{u+1}}\sum_{n\in C_N(u,a)}
 \frac{1}{\norm{n\alpha}^{s}}.\]
For  \(n\in C_N(u,a)\), Lemma \ref{lem:two-sided-norm-bound} yields
\[
    \norm{n\alpha}\ge a|D_u|\ge\frac{a}{2q_{u+1}}.
\]
Furthermore, the set \(C_N(u,a)\) can be explicitly decomposed into a disjoint union of two truncated sets:
\begin{equation*}
 C_N(u,a)=
 A_N(\underbrace{0,\ldots,0}_{u \text{ times}},a+1)
 \sqcup A_N(\underbrace{0,\ldots,0}_{u-1 \text{ times}},1,a_{u+1}-a),
\end{equation*}
where \(A_N(\cdot)\)  defined in \eqref{cylinder}.
Applying Lemma \ref{lem:card-bound-Ostrowski-A_N},  we obtain 
\begin{equation*}
 \sharp C_N(u,a)\le2\left(\frac{3N}{q_{u+1}}+1\right)
 \le\frac{8N}{q_{u+1}},\qquad u\le K-1.
\end{equation*}

Hence, for $0<s<1$, we have
\begin{align*}
  \sum_{\substack{n \in \Omega_{N,1} \\ \norm{n\alpha} < |D_2|}} \frac{1}{\norm{n\alpha}^{s}}
 \ll N\sum_{u=2}^{K-1}\frac{q_{u+1}^{s}}{q_{u+1}}
 \sum_{a=1}^{a_{u+1}}\frac{1}{a^s}
 \ll N\sum_{u=2}^{K-1}\frac{a_{u+1}^{1-s}}{q_{u+1}^{1-s}}
 \ll N.
\end{align*}
For \(s>1\), we have
\begin{align*}
   \sum_{\substack{n \in \Omega_{N,1} \\ \norm{n\alpha} < |D_2|}} \frac{1}{\norm{n\alpha}^{s}}
 \ll N\sum_{u=2}^{K-1}\frac{q_{u+1}^{s}}{q_{u+1}}
 \sum_{a=1}^{a_{u+1}}\frac{1}{a^s}
 \ll N\sum_{u=2}^{K-1}q_{u+1}^{s-1}
 \ll Nq_K^{s-1}.
\end{align*}

    \end{proof}
\begin{lem}\label{R_2u}
    For any \(s>0\), we have
    \[
         \sum_{\substack{1\le n\le N \\\norm{n\alpha} \ge |D_2|}} \frac{1}{\norm{n\alpha}^{s}} \le 2q_3^s N.
    \]
\end{lem}
\begin{proof}
A direct calculation yields
\begin{equation*}
    \sum_{\substack{1\le n\le N \\ \norm{n\alpha}\ge|D_2|}} \frac{1}{\norm{n\alpha}^{s}}
    \le \frac{N}{|D_2|^s} \le 2q_3^s N.
\end{equation*}
    
\end{proof}

Combining Lemmas \ref{R_1u} and \ref{R_2u}, we obtain
\begin{equation}\label{u}
    \sum_{n \in \Omega_{N,1}} \frac{1}{\norm{n\alpha}^{s}}\ll  
    \begin{cases}
        N+q_3^{s}N, & \text{if } 0<s<1, \\[1.4em]
        Nq_K^{s-1}+q_3^sN, & \text{if } s>1.
    \end{cases}
\end{equation}

\noindent\underline{\textbf{Lower bounds for \(\sum_{n \in \Omega_{N,1}}\norm{n\alpha}^{-s}\).}}
For $m,a\in\N$, recall that
\begin{equation*}
 A_N(\underbrace{0,\ldots,0}_{m \text{ times}},a)=\left\{1\le n\le N:
 \begin{array}{l}
 \text{the Ostrowski expansion of }n\text{ satisfies}\\
 c_{k+1}=0\text{ for all }k<m,\quad c_{m+1}=a
 \end{array}\right\}.
\end{equation*}
To simplify notation, we write \(A(0^m, a)\) for \(A_N(0,\ldots,0,a)\). Since $1\le n\le N<q_{K+1}$, we have $m\le K$. Note that
if \(m=K\), then any \(n\in A(0^m, a)\) must be of the form \(n=aq_K\), which yields \(n\equiv0\pmod{q_K}\). Thus, we restrict our attention to \(m \le K-1\).

Furthermore, it is easily checked that if \(n \in A(0^m, a))\) satisfies one of the following conditions:
\begin{itemize}
    \item \(1 \le a \le a_{m+1}\) and \(m \le K - 2\), or
    \item \(2 \le a \le a_{m+1}\) and \(m = K - 1\),
\end{itemize}
then \(n \in \Omega_{N,1}\).  

Hence, for $0<s<1$, we obtain 
\begin{align}
 \sum_{n\in\Omega_{N,1}}
 \frac{1}{\norm{n\alpha}^{s}}&\ge\sum_{m=2}^{K-2}\sum_{a=1}^{a_{m+1}}
 \sum_{n\in A(0^m, a)}\frac{1}{\norm{n\alpha}^{s}}
 +\sum_{a=2}^{a_K}\sum_{n\in A(0^{K-1}, a)}\frac{1}{\norm{n\alpha}^{s}}
 \notag\\
 &\gg\sum_{m=2}^{K-2}\sum_{a=1}^{a_{m+1}}
 \frac{N}{q_{m+1}}\frac{1}{a^s|D_m|^s}
 +\sum_{a=2}^{a_K}\frac{N}{q_K}\frac{1}{a^s|D_{K-1}|^s}\notag\\
 &\gg N\left(\sum_{m=2}^{K-2}\sum_{a=1}^{a_{m+1}}
 \frac{1}{q_{m+1}^{1-s}a^s}
 +\sum_{a=2}^{a_K}\frac{1}{q_K^{1-s}a^s}\right)\notag\\
 &=N\left(\sum_{m=2}^{K-1}\sum_{a=1}^{a_{m+1}}
 \frac{1}{q_{m+1}^{1-s}a^s}-\frac{1}{q_K^{1-s}}\right)\notag\\
 &\gg N\sum_{m=2}^{K-1}\frac{a_{m+1}^{1-s}}{q_{m+1}^{1-s}}
 -\frac{N}{q_K^{1-s}}\notag\\
 &\gg N\left(\sum_{m=2}^{K-1}\frac{1}{q_m^{1-s}}
 -\frac{1}{q_K^{1-s}}\right)\notag\\
 &\gg N\sum_{m=2}^{K}\frac{1}{q_m^{1-s}}
 \gg N,\label{eq:regular-lower-subcritical-chain}
\end{align}
where the second inequality follows from Lemmas \ref{lem:two-sided-norm-bound} and \ref{lem:card-bound-Ostrowski-A_N}.

For $s>1$, we have
\begin{align}
 \sum_{n\in\Omega_{N,1}}
 \frac{1}{\norm{n\alpha}^{s}}&\ge\sum_{m=2}^{K-2}\sum_{a=1}^{a_{m+1}}
 \sum_{n\in A(0^m, a)}\frac{1}{\norm{n\alpha}^{s}}
 +\sum_{a=2}^{a_K}\sum_{n\in A(0^{K-1}, a)}\frac{1}{\norm{n\alpha}^{s}}
 \notag\\
 &\gg\sum_{a=2}^{a_K}\frac{N}{q_K}\frac{1}{a^s|D_{K-1}|^s}\asymp
 Nq_K^{s-1}.\label{eq:regular-lower-supercritical-chain}
\end{align}

Combining the upper bounds \eqref{u}  with the lower bounds \eqref{eq:regular-lower-subcritical-chain} and \eqref{eq:regular-lower-supercritical-chain}, we obtain
\begin{equation}\label{r}
    \sum_{n \in \Omega_{N,1}} \frac{1}{\norm{n\alpha}^{s}}\asymp
    \begin{cases}
        N, & \text{if } 0<s<1, \\[1.4em]
        Nq_K^{s-1}, & \text{if } s>1,
    \end{cases}
\end{equation}
which completes the proof for the regular classes.

\subsection{Estimates for the exceptional classes}\label{subsec:exceptional}

We now treat the exceptional  classes $n\equiv0,q_{K-1}\pmod{q_K}$. We begin with the class \(n\equiv0\pmod{q_K}\).

First, note that \(N\ge q_3\) implies \(K\ge3\). Let \(1\le n\le N\) satisfy \(n\equiv0\pmod{q_K}\). Then
\[
    n = aq_K \quad \text{with} \quad 1 \le a \le \frac{N}{q_K}.
\]
Further, since \(N < q_{K+1} = a_{K+1}q_K + q_{K-1}\), we have
\[
    a \le \frac{N}{q_K} < a_{K+1} + \frac{q_{K-1}}{q_K} < a_{K+1} + 1,
\]
which implies \(a \le a_{K+1}\). Thus, \(n = aq_K\) is the Ostrowski expansion of \(n\).

Since the only non-zero Ostrowski coefficient is \(c_{K+1}=a\), it follows from Lemma \ref{lem:homogeneous-norm-estimate} that
\[
    \norm{n\alpha} = a|D_K|\asymp\frac{a}{q_{K+1}}.
\]
For \(0<s<1\), this gives
\begin{align}
    \sum_{\substack{1\le n\le N \\ n\equiv0\pmod{q_K}}}
    \frac{1}{\norm{n\alpha}^{s}}
    \asymp q_{K+1}^s\sum_{1\le a\le N/q_K}\frac{1}{a^s}
    \asymp \frac{N^{1-s}q_{K+1}^{s}}{q_K^{1-s}}.
    \label{eq:multiple-sum-subcritical}
\end{align}
For \(s>1\), we similarly obtain
\begin{align}
    \sum_{\substack{1\le n\le N \\ n\equiv0\pmod{q_K}}}
    \frac{1}{\norm{n\alpha}^{s}}
    \asymp q_{K+1}^s\sum_{1\le a\le N/q_K}\frac{1}{a^s}
    \asymp q_{K+1}^{s}.
    \label{eq:multiple-sum-supercritical}
\end{align}

For the shifted residue class \(n \equiv q_{K-1} \pmod{q_K}\), we write \(n\) as
\begin{equation*}
    n=q_{K-1}+(a_{K+1}-b)q_K \quad \text{with} \quad
    a_{K+1}-\frac{N-q_{K-1}}{q_K}\le b\le a_{K+1}
\end{equation*}
for some integer \(b\). Since \(N<q_{K+1}=a_{K+1}q_K+q_{K-1}\), we have
\begin{equation*}
    a_{K+1}-\frac{N-q_{K-1}}{q_K}>0,
    \qquad\text{and hence } b\ge1.
\end{equation*}
Since \(K\ge3\), we have \(K-1\ge2\). Applying Lemma \ref{lem:homogeneous-norm-estimate} yields
\begin{align*}
    \norm{n\alpha}
    &=|D_{K-1}|-(a_{K+1}-b)|D_K| \\
    &=|D_{K+1}|+b|D_K|\asymp b|D_K|\asymp\frac{b}{q_{K+1}}.
\end{align*}
Thus, for \(0<s<1\),
\begin{align}
    \sum_{\substack{1\le n\le N \\ n\equiv q_{K-1}\pmod{q_K}}}
    \frac{1}{\norm{n\alpha}^{s}}
    &\ll\sum_{\substack{b\in\mathbb{Z}\\a_{K+1}-(N-q_{K-1})/q_K\le b\le a_{K+1}}}
    \frac{q_{K+1}^{s}}{b^s}\notag\\
    &\ll q_{K+1}^{s}\sum_{1\le b\le N/q_K}\frac{1}{b^s}
    \ll\frac{q_{K+1}^{s}N^{1-s}}{q_K^{1-s}},
    \label{eq:shifted-sum-subcritical}
\end{align}
and for \(s>1\), we similarly obtain
\begin{equation}\label{eq:shifted-sum-supercritical}
    \sum_{\substack{1\le n\le N \\ n\equiv q_{K-1}\pmod{q_K}}}
    \frac{1}{\norm{n\alpha}^{s}}
    \ll\sum_{\substack{b\in\mathbb{Z}\\a_{K+1}-(N-q_{K-1})/q_K\le b\le a_{K+1}}}
    \frac{q_{K+1}^{s}}{b^s}\ll q_{K+1}^{s}.
\end{equation}

Combining \eqref{eq:multiple-sum-subcritical} and \eqref{eq:shifted-sum-subcritical} for \(0<s<1\), and \eqref{eq:multiple-sum-supercritical} and \eqref{eq:shifted-sum-supercritical} for \(s>1\), we conclude that
\begin{equation}\label{e}
    \sum_{n\in\Omega_{N,2}}
    \frac{1}{\norm{n\alpha}^{s}}
    \asymp
    \begin{cases}
        \displaystyle \frac{N^{1-s}q_{K+1}^{s}}{q_K^{1-s}}, & \text{if } 0<s<1, \\[1.6em]
        \displaystyle q_{K+1}^{s}, & \text{if } s>1.
    \end{cases}
\end{equation}

Finally, together with \eqref{r} and \eqref{e}, this completes the proofs of Propositions \ref{pro:s<1} and \ref{pro:s>1}. In what follows,
we obtain the following asymptotic formula for the full reciprocal sum:
\begin{equation}\label{total}
    \sum_{1\le n\leq N} \frac{1}{\norm{n\alpha}^s} \asymp 
    \begin{cases}
        N + \displaystyle\frac{N^{1-s}q_{K+1}^{s}}{q_K^{1-s}}, & \text{if } 0<s<1, \\[1.8em]
        N q_K^{s-1} + q_{K+1}^{s}, & \text{if } s>1.
    \end{cases}
\end{equation}

\section{Proof of Theorem \ref{thm:metric-subcritical}}
For clarity, the proof is divided into the following two parts:
\begin{itemize}
    \item the convergence part for \(0<s<1\),
    \item the divergence part for \(0<s<1\).
\end{itemize}

\noindent\underline{Convergence part \((0<s<1)\).} Choose a sufficiently small \(\varepsilon_0>0\) such that \(w(\alpha)+\varepsilon_0<1/s\). Then, by Lemma \ref{lem:omega-def}, for all sufficiently large \(K\), we have
\begin{equation*}
    q_{K+1} \le q_K^{w(\alpha)+\varepsilon_0} \le N^{w(\alpha)+\varepsilon_0} < N^{1/s}.
\end{equation*}
Combining this with \eqref{total}, we obtain
\[
    \sum_{1\le n\leq N} \frac{1}{\norm{n\alpha}^s} \asymp N.
\]

Consequently, applying Abel's summation formula yields
\begin{align*}
    \sum_{1\le n\le N}\frac{\psi(n)}{\norm{n\alpha}^{s}}
    &= \sum_{n=1}^{N-1}\bigl(\psi(n)-\psi(n+1)\bigr) \sum_{m=1}^{n}\frac{1}{\norm{m\alpha}^{s}}
    + \psi(N)\sum_{m=1}^{N}\frac{1}{\norm{m\alpha}^{s}}\\
    &\ll \sum_{n=1}^{N-1}\bigl(\psi(n)-\psi(n+1)\bigr) n 
    + \psi(N) N \\
    &= \sum_{1\le n\le N}\psi(n).
\end{align*}
If \(\sum \psi(n) < \infty\), the first Borel--Cantelli lemma yields \(\mbox{\(\mathcal{L}\bigl(W_s(\psi,\alpha,\gamma)\bigr) = 0\)}\), completing the proof of the convergence part.

\medskip
\noindent\underline{Divergence part \((0<s<1)\).}
Since \(\psi(n)\le \psi(n)/\norm{n\alpha}^{s}\), we have the inclusion \(W(\psi,\gamma)\subset W_s(\psi,\alpha,\gamma)\). By Theorem \ref{thm:inhomogeneous-Khintchine}, the divergence condition \(\sum_{n=1}^{\infty}\psi(n)=\infty\) ensures that \(W(\psi,\gamma)\) has full measure. It follows immediately that 
\(
    \mathcal{L}\bigl(W_s(\psi,\alpha,\gamma)\bigr)=1.
\)

\section{Proof of Theorem \ref{thm:metric-supercritical-initial}}

The proof proceeds in two steps:
\begin{itemize}
    \item the convergence part for \(s>1\),
    \item the divergence part for \(s>1\).
\end{itemize}

\noindent\underline{Convergence part \((s>1)\).} For \(\varepsilon>0\), combining Lemma \ref{lem:omega-def} with \(q_K\le N<q_{K+1}\), we deduce that for all sufficiently large \(K\),
\begin{equation*}
    q_{K+1}\le q_K^{w(\alpha)+\varepsilon/s} \le N^{w(\alpha)+\varepsilon/s}.
\end{equation*}
Substituting this bound into the estimate of (\ref{total}), we obtain
\begin{equation*}
 \sum_{1\le n\le N}\frac1{\norm{n\alpha}^{s}}
 \ll q_{K+1}^{s}\le N^{w(\alpha)s+\varepsilon}.
\end{equation*}
Combining this bound with Abel's summation formula, we deduce that
\begin{align*}
    \sum_{1\le n\le N}\frac{\psi(n)}{\norm{n\alpha}^{s}}
    &= \sum_{n=1}^{N-1}\bigl(\psi(n)-\psi(n+1)\bigr) \sum_{m=1}^{n}\frac{1}{\norm{m\alpha}^{s}}
    + \psi(N)\sum_{m=1}^{N}\frac{1}{\norm{m\alpha}^{s}} \notag \\
    &\ll \sum_{n=1}^{N-1}\bigl(\psi(n)-\psi(n+1)\bigr) n^{w(\alpha)s+\varepsilon}
    + \psi(N) N^{w(\alpha)s+\varepsilon} \notag \\
    &\asymp  \sum_{1\le n\le N}\psi(n)n^{w(\alpha)s+\varepsilon-1}.
\end{align*}
Therefore, if for some \(0<\varepsilon_0<1\) we have
\begin{equation*}
    \sum_{n=1}^{\infty}\psi(n)n^{w(\alpha)s+\varepsilon_0-1}<\infty,
\end{equation*}
then the first Borel--Cantelli lemma gives
\begin{equation*}
    \mathcal{L}\bigl(W_s(\psi,\alpha,\gamma)\bigr)=0.
\end{equation*}

\medskip

\noindent\underline{Divergence part \((s>1)\).}
According to Theorem \ref{thm:OpenAI}, it suffices to establish the following estimate:
\begin{equation*}
    \sum_{1\le n\le N}\frac{\phi(n)}{n} \frac{\psi(n)}{\norm{n\alpha}^{s}}
    \gg \sum_{1\le n\le N}\frac{\psi(n)n^{s-1}}{\log\log n}.
\end{equation*}

We first establish two-sided bounds for the unweighted sum.

\begin{lem}\label{lem:totient-reciprocal-two-sided}
Let \(C_{\alpha}>0\) be a constant depending only on \(\alpha\). We have
\begin{equation*}
    \frac{\phi(q_K)}{q_K}\,q_{K+1}^{s}
    \le\sum_{1\le n\le N}\frac{\phi(n)}{n\norm{n\alpha}^{s}}
    \le C_{\alpha}q_{K+1}^{s}.
\end{equation*}
\end{lem}
\begin{proof}
The lower bound follows by choosing \(n=q_K\). The upper bound is a direct consequence of \eqref{total}.
\end{proof}

\begin{lem}\label{lem:totient-reciprocal-lower}
    For \(N \ge q_3\), we have
    \begin{equation*}
        \sum_{1\le n\le N}\frac{\phi(n)}{n\norm{n\alpha}^{s}} \gg\frac{N^s}{\log\log N}.
    \end{equation*}
\end{lem}

\begin{proof}
Combining the lower bound in Lemma \ref{lem:totient-reciprocal-two-sided} with the classical estimate \cite[pp.~216--219]{Landau1909}: \(\phi(n)/n \gg 1/\log\log n\), we deduce that
\begin{equation*}
    \sum_{1\le n\le N}\frac{\phi(n)}{n\norm{n\alpha}^{s}}
    \ge\frac{\phi(q_K)}{q_K}q_{K+1}^{s}
    \gg\frac{N^s}{\log\log q_K}
    \gg\frac{N^s}{\log\log N}.
\end{equation*}
\end{proof}

Finally, applying Abel's summation formula and Lemma \ref{lem:totient-reciprocal-lower}, we obtain the desired weighted lower bound.

\begin{lem}\label{lem:totient-weighted-comparison}
\begin{equation*}
    \sum_{1\le n\le N}\frac{\phi(n)}{n} \frac{\psi(n)}{\norm{n\alpha}^{s}}
    \gg\sum_{1\le n\le N}\frac{\psi(n)n^{s-1}}{\log\log n}.
\end{equation*}
\end{lem}
\begin{proof}
Since \(\psi\) is nonincreasing, Abel's summation formula yields
\begin{align*}
    \sum_{1\le n\le N} \frac{\phi(n)}{n} \frac{\psi(n)}{\norm{n\alpha}^{s}}
    &= \sum_{n=1}^{N-1}\bigl(\psi(n)-\psi(n+1)\bigr) \sum_{m=1}^{n}\frac{\phi(m)}{m\norm{m\alpha}^{s}} + \psi(N)\sum_{m=1}^{N}\frac{\phi(m)}{m\norm{m\alpha}^{s}} \\
    &\gg \sum_{n=1}^{N-1}\bigl(\psi(n)-\psi(n+1)\bigr) \frac{n^s}{\log\log n} + \psi(N) \frac{N^s}{\log\log N} \\
    &\asymp \sum_{1\le n\le N}\psi(n)\frac{n^{s-1}}{\log\log n}.
\end{align*}
\end{proof}

Combining Theorem \ref{thm:OpenAI} and Lemma \ref{lem:totient-weighted-comparison}, we complete the proof of the divergence part.


\begin{thebibliography}{99}



\bibitem{ANL21}
F.~Adiceam, E.~Nesharim and F.~Lunnon.
\textit{On the \(t\)-adic Littlewood conjecture.}
\newblock \emph{Duke Math. J.}
\textbf{170} (2021), no.~10, 2371--2419.







\bibitem{BHV13}
V.~Beresnevich, A.~Haynes and S.~Velani.
\textit{Multiplicative zero-one laws and metric number theory.}
\newblock \emph{Acta Arith.}
\textbf{160} (2013), no.~2, 101--114.




\bibitem{Beresnevich-Haynes-Velani2020}
V.~Beresnevich, A.~Haynes, and S.~Velani.
\textit{Sums of reciprocals of fractional parts and multiplicative Diophantine approximation.}
\newblock \emph{Mem. Amer. Math. Soc.}
\textbf{263} (2020), no.~1276.






\bibitem{CS55}
J.~Cassels and H.~Swinnerton-Dyer.
\textit{On the product of three homogeneous linear forms and indefinite ternary quadratic forms.}
\newblock \emph{Philos. Trans. Roy. Soc. London Ser. A}
\textbf{248} (1955), 73--96.



\bibitem{Chow2018}
S.~Chow.
\textit{Bohr sets and multiplicative Diophantine approximation.}
\newblock \emph{Duke Math. J.}
\textbf{167} (2018), no.~9, 1623--1642.



\bibitem{CHPR25}
S.~Chow, M.~Hauke, A.~Pollington and F.~A.~Ram\'irez.
\textit{General Duffin--Schaeffer-type counterexamples in Diophantine approximation.}
\newblock \emph{Mathematika}, to appear, preprint 2025, arXiv:2504.16565.


\bibitem{CT19}
S.~Chow and N.~Technau.
\textit{Higher-rank Bohr sets and multiplicative Diophantine approximation.}
\newblock \emph{Compos. Math.}
\textbf{155} (2019), no.~11, 2214--2233.


\bibitem{Chow-Zafeiropoulos2021}
S.~Chow and A.~Zafeiropoulos.
\textit{Fully-inhomogeneous multiplicative Diophantine approximation of badly approximable numbers.}
\newblock \emph{Int. Math. Res. Not. IMRN}
(2022), no.~15, 11637--11670.


\bibitem{Chow-Technau2024}
S.~Chow and N.~Technau.
\textit{Littlewood and Duffin--Schaeffer-Type Problems in Diophantine Approximation.}
\newblock \emph{Mem. Amer. Math. Soc.}
\textbf{296} (2024), no.~1475.

\bibitem{Chow-Yang2024}
S.~Chow and L.~Yang.
\textit{Effective equidistribution for multiplicative Diophantine approximation on lines.}
\newblock \emph{Invent. Math.}
\textbf{235} (2024), no.~3, 973--1007.

\bibitem{Duffin-Schaeffer1941}
R.~Duffin and A.~Schaeffer.
\textit{Khintchine's problem in metric Diophantine approximation.}
\newblock \emph{Duke Math. J.}
\textbf{8} (1941), 243--255.



\bibitem{EKL06}
M.~Einsiedler, A.~Katok and E.~Lindenstrauss.
\textit{Invariant measures and the set of exceptions to Littlewood's conjecture.}
\newblock \emph{Ann. of Math. (2)}
\textbf{164} (2006), 513--560.


\bibitem{E70}
P.~Erd\H{o}s,
\textit{On the distribution of the convergents of almost all real numbers},
\newblock \emph{J. Number Theory}
\textbf{2} (1970), 425--441.


\bibitem{FH24}
L.~Fr\"uhwirth and M.~Hauke.
\textit{The Duffin-Schaeffer Conjecture for multiplicative Diophantine approximation.}
\newblock (2024), arXiv: 2403.11257.



\bibitem{G61}
P.~Gallagher,
\textit{Approximation by reduced fractions.}
\newblock \emph{J. Math. Soc. Japan}
\textbf{13} (1961), 342--345.

\bibitem{Gallagher1962}
P.~Gallagher.
\textit{Metric simultaneous Diophantine approximation.}
\newblock \emph{J. London Math. Soc.}
\textbf{37} (1962), 387--390.


\bibitem{Gorodnik-Vishe2018}
A.~Gorodnik and P.~Vishe.
\textit{Diophantine approximation for products of linear maps---logarithmic improvements.}
\newblock \emph{Trans. Amer. Math. Soc.}
\textbf{370} (2018), 487--507.




\bibitem{Harman1998}
G.~Harman.
\textit{Metric Number Theory.}
\newblock London Mathematical Society Monographs. New Series, vol.~18,
Clarendon Press, Oxford University Press, New York, 1998.


\bibitem{HMP26}
M.~Hauke-Treuer, J.~Maynard, and A.~Pollington.
\textit{Counterexamples to the inhomogeneous Duffin--Schaeffer conjecture.}
\newblock arXiv:2609.30921v1, 2026. 




\bibitem{HL26}
Y.~He and L.~Liao.
\textit{Counterexamples to the inhomogeneous Duffin--Schaeffer conjecture for a residual set of shifts.}
\newblock arXiv:2609.30870v1, 2026. 


\bibitem{Khintchine1924}
A.~Khintchine.
\textit{Einige S\"atze \"uber Kettenbr\"uche, mit Anwendungen auf die Theorie der Diophantischen Approximationen.}
\newblock \emph{Math. Ann.}
\textbf{92} (1924), 115--125.


\bibitem{Koukoulopoulos-Maynard2020}
D.~Koukoulopoulos and J.~Maynard.
\textit{On the Duffin--Schaeffer conjecture.}
\newblock \emph{Ann. of Math. (2)}
\textbf{192} (2020), 251--307.





\bibitem{Landau1909}
E.~Landau.
\textit{Handbuch der Lehre von der Verteilung der Primzahlen.}
\newblock Vol.~1, Teubner, Leipzig, 1909.



\bibitem{openai2026}
OpenAI.
\textit{The Weak Inhomogeneous Duffin--Schaeffer Conjecture.}
\newblock OpenAI Math Release preprint,
OAI: The-weak-inhomogeneous-Duffin-Schaeffer-conjecture-September-25-2026, 2026.


\bibitem{Ostrowski1922}
A.~Ostrowski.
\textit{Bemerkungen zur Theorie der Diophantischen Approximationen.}
\newblock Abh. Math. Semin. Univ. Hambg., 1:77--98, 1922.


\bibitem{P61}
L.~G.~Peck.
\textit{Simultaneous rational approximations to algebraic numbers.}
\newblock \emph{Bull. Amer. Math. Soc.}
\textbf{67} (1961), 197--201.



\bibitem{PV90}
A.~D.~Pollington and R.~C.~Vaughan,
\textit{The \(k\)-dimensional Duffin and Schaeffer conjecture}.
\newblock Mathematika
\textbf{37} (1990), 190--200.


\bibitem{Pollington-Velani2000}
A.~D.~Pollington and S.~L.~Velani.
\textit{On a problem in simultaneous Diophantine approximation: Littlewood's conjecture.}
\newblock \emph{Acta Math.}
\textbf{185} (2000), 287--306.



\bibitem{R17}
F.~A.~Ram\'irez.
\textit{Counterexamples, covering systems, and zero-one laws for inhomogeneous approximation.}
\newblock \emph{Int. J. Number Theory}
\textbf{13} (2017), 633--654.


\bibitem{RS1992}
A.~M.~Rockett and P.~Sz\"usz.
\textit{Continued fractions.}
\newblock World Scientific Publishing Co., Inc., River Edge, NJ, 1992.


\bibitem{Shapira2011}
U.~Shapira.
\textit{A solution to a problem of Cassels and Diophantine properties of cubic numbers.}
\newblock \emph{Ann. of Math. (2)}
\textbf{173} (2011), 543--557.


\bibitem{Szusz1958}
P.~Sz\"usz.
\textit{\"Uber die metrische Theorie der Diophantischen Approximation.}
\newblock \emph{Acta Math. Hungar.}
\textbf{9} (1958), 177--193.


\bibitem{S60}
W.~M.~Schmidt.
\textit{On normal numbers.}
\newblock \emph{Pacific J. Math.}
\textbf{10} (1960), 661--672.



\bibitem{V78}
J.~D.~Vaaler,
\textit{On the metric theory of Diophantine approximation},
\newblock \emph{Pacific J. Math.}
\textbf{76} (1978), 527--539.

\bibitem{Yu2019}
H.~Yu.
\textit{A Fourier-analytic approach to inhomogeneous Diophantine approximation.}
\newblock \emph{Acta Arith.}
\textbf{190} (2019), 263--292.



\bibitem{Yu2021}
H.~Yu.
\textit{On the metric theory of inhomogeneous Diophantine approximation: an Erd\H{o}s-Vaaler type result.}
\newblock \emph{J. Number Theory}
\textbf{224} (2021), 243--273.



\bibitem{Yu2023}
H.~Yu.
\textit{On the metric theory of multiplicative Diophantine approximation.}
\newblock \emph{J. Anal. Math.}
\textbf{149} (2023), no.~2, 763--800.
\end{thebibliography}
\end{document}